\documentclass[11pt]{article}

\usepackage{amsfonts,amsmath,amsthm}
\usepackage{amssymb}
\usepackage{latexsym}

\usepackage[english]{babel}

\newtheoremstyle{RemarkExampleStyle1}{3pt}{3pt}{}{}{}{.}{6pt}{}

\theoremstyle{plain}
\newtheorem{thm}{Theorem}[section]

\newtheorem{prop}[thm]{Proposition}
\newtheorem{lemma}[thm]{Lemma}
\newtheorem{cor}[thm]{Corollary}
\newtheorem{definition}[thm]{Definition}
\newtheorem*{clm}{Claim}

\theoremstyle{RemarkExampleStyle1}
\newtheorem{remark}[thm]{R e m a r k }
\newtheorem*{example}{E x a m p l e}

\newcounter{abc}   
\newcounter{iiiii} 

\newenvironment{aequivalenz}
{\setcounter{iiiii}{0}
\begin{list}%
{{\rm (\roman{iiiii})}}
{\usecounter{iiiii}
\parsep=0pt plus 1pt
\topsep=1pt plus 2pt minus 1pt
\itemsep=1pt plus 2pt minus 1pt
\leftmargin=3\baselineskip \labelsep=.6\baselineskip
\labelwidth=2.4\baselineskip
\rightmargin 0pt}
}
{\end{list}}

\newenvironment{statements}
{\setcounter{abc}{0}
\begin{list}
{{\rm (\alph{abc})}}
{\usecounter{abc}
\parsep=0pt plus 1pt
\topsep=1pt plus 2pt minus 1pt
\itemsep=1pt plus 2pt minus 1pt
\leftmargin=3\baselineskip \labelsep=.6\baselineskip
\labelwidth=2.4\baselineskip
\rightmargin 0pt}
}
{\end{list}}

\numberwithin{equation}{section}

\begin{document}
\sloppy
\title{\sc Daugavet centers}

\author{T. Bosenko and V. Kadets\footnote{Research of the second named author was conducted
 during his stay in the University of Granada and was
 supported by Junta de Andaluc\'{\i}a grant P06-FQM-01438.}}

\date{}




\maketitle

\thispagestyle{empty}

\vspace{-25pt}

\begin{center}
\small 

{\it Department of Mechanics and Mathematics,
V.N.~Karazin Kharkov National University,
 4~Svobody~Sq., Kharkov, 61077, Ukraine}

\bigskip
E-mail: t.bosenko@mail.ru

E-mail: vova1kadets@yahoo.com
\end{center}

\begin{abstract}

An operator $G {:}\allowbreak\  X \to Y$ is said to be a Daugavet center if $\|G + T\|
= \|G\| + \|T\|$ for every rank-$1$ operator $T {:}\allowbreak\  X \to Y$.
The main result of the paper is: if $G {:}\allowbreak\  X \to Y$ is a Daugavet center,
$Y$ is a subspace of a Banach space \, $E$, and $J: Y \to E$ is the natural
embedding operator, then $E$ can be equivalently renormed in such
a way, that $J \circ G : X \to E$ is also a Daugavet center. This result was
previously known for particular case $X=Y$, $G=\mathrm{Id}$ and only in
separable spaces. The proof of our generalization is based on an
idea completely different from the original one. We give also some
geometric characterizations of Daugavet centers, present a number of
examples, and generalize (mostly in straightforward manner) to Daugavet centers
some results known previously for spaces with the Daugavet property.

\end{abstract}

\smallskip



\smallskip


\section{Introduction}

A Banach space $X$ is said to have the Daugavet property
if all the operators $T {:}\allowbreak\  X \to X$ of rank-$1$ satisfy
the Daugavet equation
 \begin{equation} \label{0-eq1}
  \|\mathrm{Id}+T\|=1+\|T\|
 \end{equation}
Several classical spaces have the Daugavet property: $C(K)$ where
K is perfect \cite{Daug}, $L_{1}(\mu)$ where $\mu$ has no atoms
\cite{Loz}, and certain functional algebras such as the disk
algebra $A(\mathbb D)$ or the algebra of bounded analytic
functions $H^\infty$ (\cite{Woj92}, \cite{Dirk10}).

 Geometric and linear-topological properties of such spaces
 were studied intensively during the last two decades (see the
 survey paper \cite{dirk-irbull} and most recent developments in
 \cite{KW}, \cite{KadKalWer} and \cite{KadShepW}). In particular,
 if $X$ is a space with the Daugavet property then every weakly compact operator,
 even every strong Radon-Nikod\'ym operator on $X$, and every operator on $X$
not fixing a copy of $\ell_{1}$ fulfill (\ref{0-eq1}) as well
(\cite{KadSSW}, \cite{Shv1}). These spaces contain subspaces
isomorphic to $\ell_{1}$, cannot have the Radon-Nikod\'ym
property, never have an unconditional basis and even never embed
into a space having an unconditional basis. The key to the later
embedding property is the following theorem:

\begin{thm}\label{theo-intr-1} \cite[Theorem 2.5]{KadSSW}
Let $X$ be a subspace of a separable Banach space $Y$, $J: X \to Y$ be the
natural embedding operator, and suppose $X$ has the Daugavet property. Then $Y$
can be renormed so that the new norm coincides with the original
one on $X$ and in the new norm $\|J+T\|=1+\|T\|$ for every
rank-$1$ operator $T {:}\allowbreak\  X \to Y$.
\end{thm}

The aim of our paper is to take off the separability condition in
the above theorem. On this way we introduce and study the
following concept:

\begin{definition} \label{defDC}
Let $X$ and $Y$ be Banach spaces. A linear continuous operator $G
{:}\allowbreak\  X \to Y$ is said to be a Daugavet center if the norm identity
\begin{equation}  \label{eqDC}
\|G + T\| = \|G\| + \|T\|
\end{equation}
is fulfilled for every rank-$1$ operator $T {:}\allowbreak\  X \to Y$.
\end{definition}

Our main result is more general than just the non-separable
version of Theorem \ref{theo-intr-1}. Namely we prove the
following:

\begin{thm}\label{theo-intr-2}
If $G {:}\allowbreak\  X \to Y$ is a Daugavet center, $Y$ is a subspace of a Banach space \ $E$,
and $J: Y \to E$ is the natural embedding operator, then $E$ can
be equivalently renormed in such a way, that the new norm
coincides with the original one on $Y$, and $J \circ G : X \to E$
is also a Daugavet center.
\end{thm}

Let us explain the structure of the paper. In Section~\ref{sec2}
of this paper we collect some straightforward generalizations to
Daugavet centers of the properties known for $\mathrm{Id}$ in the spaces with the
Daugavet property. We study also properties of the unit ball images
under Daugavet centers. In Section~\ref{sec3} we
give some examples of Daugavet centers quite different from the identity
operator or isometric embedding, which were known before. Finally, Section~\ref{main} is devoted to the proof of the main result.

In this paper we deal with real Banach spaces. We use the letters $X, Y, E$ to denote Banach spaces and their
subspaces. $L(X,Y)$ stands for the space of all linear bounded
operators, acting from $X$ to $Y$. $B_X$ denotes the closed unit ball of
a Banach space $X$ and $S_X$ denotes its unit sphere. For a bounded closed convex set $A \subset X$ and for $x^* \in X^*$ we denote
$$S(A, x^*, \varepsilon )=\{x \in A{:}\allowbreak\  x^*(x) \geq \sup x^*(A)-\varepsilon\}$$
the slice of $A$, generated by $x^*$.
We use the notation
$$
S(x^*, \varepsilon)= \{x\in B_X{:}\allowbreak\  x^*(x)\ge 1-\varepsilon \}
$$
for the slice of $B_X$ determined by a functional $x^*\in S_{X^*}$
and $\varepsilon>0$.

We say that an element $x \in A$  is a denting point of the set
$A$ if for every $\varepsilon>0$ there is a slice of $A$ which contains
$x$ and has diameter smaller than $\varepsilon$.

A set $A$ is said to have the Radon-Nikod\'ym property if every
closed convex subset $B \subset A$ is the closed convex hull of
its denting points.

Operator $T \in L(X, Y)$ is said to be a strong Radon-Nikod\'ym
operator if the closure of $T(B_X)$ has the Radon-Nikod\'ym
property.

\smallskip

\section{Basic properties of Daugavet centers} \label{sec2}

Definition \ref{defDC} implies the equality $\|aG + bT\| = a\|G\|
+b\|T\|$ for every $a, b \geq 0$. This means that a non-zero
operator $G$ is a Daugavet center if and only if  $G/ \|G\|$ is. Therefore below we
mostly consider the case $\|G\|=1$, and by the same reason when it
is convenient, we require $\|T\|=1$.

\begin{thm} \label{theo2.1}
For an operator $G \in L(X,Y)$ with $\|G\|=1$ the following
assertions are equivalent:
\begin{aequivalenz}
\item $G$ is a Daugavet center.

\item For every $y_0\in S_Y$ and every slice $S(x^*_0,\varepsilon_0)$ of
$B_X$ there is another slice $S(x^*_1, \varepsilon_1)\subset
S(x^*_0,\varepsilon_0)$ such that for every $x \in S(x^*_1, \varepsilon_1)$ the inequality $\|Gx+y_0\|>2-\varepsilon_0$ holds.

\item For every $y_0\in S_Y$, $x_0^* \in S_{X^*}$ and $\varepsilon >0$
there is $x \in B_X$ such that $x_0^*(x) \geq 1-\varepsilon$ and
$\|Gx+y_0\|>2-\varepsilon$.

\item For every $x_0^*\in S_{X^*}$ and every $weak^*$ slice
$S(B_{Y^*}, y_0, \varepsilon_0)$ (where $y_0 \in S_Y \subset
S_{Y^{**}}$) there is another $weak^*$ slice $S(B_{Y^*}, y_1, \varepsilon_1)
\subset S(B_{Y^*}, y_0, \varepsilon_0)$ such that for every $y^* \in
S(B_{Y^*}, y_1, \varepsilon_1)$ the inequality $\|G^*y^*+x_0^*\|>2-\varepsilon_0$ holds.

\item For every $x_0^*\in S_{X^*}$ and every $weak^*$ slice
$S(B_{Y^*}, y_0, \varepsilon_0)$ there is $y^* \in S(B_{Y^*}, y_0, \varepsilon_0)$
which satisfies the inequality $\|G^*y^*+x_0^*\|>2-\varepsilon_0$.

\end{aequivalenz}
\end{thm}

\noindent P r o o f. 
$(i)\Rightarrow(ii)$ Define $T$ by $Tx=x_0^*(x)y_0$. Then
$\|G^*+T^*\|=\|G+T\|=2$, so there exists a functional $y^* \in
S_{Y^*}$ such that $\| G^*y^*+T^*y^*\| \geq 2-\varepsilon_0$ and
$y^*(y_0) \geq 0$. Put
$$
 x_1^* = \frac{G^*y^*+T^*y^*}{\|G^*y^*+T^*y^*\|},
\quad \varepsilon_1 = 1- \frac{2-\varepsilon_0}{\|G^*y^*+T^*y^*\|}.
$$
Then for every $x \in S(x_1^*, \varepsilon_1)$ we have
$$
\left\langle (G^*+T^*)y^*, x \right\rangle \geq
(1-\varepsilon_1)\|G^*y^*+T^*y^*\|=2-\varepsilon_0;
$$
hence
$$
1+x_0^*(x) \geq y^*(Gx)+y^*(y_0)x_0^*(x) \geq 2-\varepsilon_0,
$$
which implies that $x_0^*(x) \geq 1-\varepsilon_0$, i.e., $x \in S(x_0^*,
\varepsilon_0)$, and
$$
2-\varepsilon_0 \leq y^*(Gx)+y^*(y_0)=y^*(Gx+y_0) \leq \|Gx+y_0\|.
$$

The implication $(ii)\Rightarrow(iii)$ is evident. Let us prove
$(iii)\Rightarrow(i)$. If $T$ is a rank-1 operator and $\|T\|=1$,
then $T$ can be represented as $Tx=x_0^*(x)y_0$ with $y_0\in S_Y$,
$x_0^* \in S_{X^*}$. Fix an $\varepsilon >0$ and let $x \in B_X$ be the
corresponding element from (iii). Then
\begin{eqnarray*}
2-\varepsilon &\leq~& \|Gx+y_0\| \leq \|Gx+x_0^*(x)y_0\|+\|(1-x_0^*(x))y_0\| \\
 &\leq&
\|(G+T)x\|+\varepsilon \|y_0\| \leq \|G+T\|+\varepsilon.
\end{eqnarray*}

So we have proved the equivalence
$(i)\Leftrightarrow(ii)\Leftrightarrow(iii)$. The remaining
equivalence $(i)\Leftrightarrow(iv)\Leftrightarrow(v)$ can be
proved the same way using dual operators $T^*$ and $G^*$ instead
of $T$ and $G$.
\bigskip

For a bounded subset $A \subset Y$ denote $r_y(A) = \sup\{\|y -
a\|{:}\allowbreak\  a \in A \}$.
 \begin{definition}\label{def2}
A bounded subset $A \subset Y$ is said to be a quasi-ball
if for every $y \in Y$
\begin{equation} \label{eq2}
r_y(A) = \|y\| + r_0(A).
\end{equation}
\end{definition}

\begin{definition}\label{def3}
A bounded subset $A \subset Y$ is called anti-dentable if
for every $y \in Y$ and for every $r \in [0, r_y(A))$
\begin{equation} \label{eq3}
  \mathop{\rm \overline{conv}}\nolimits\left( A \setminus B_Y(y,r)\right) \supset A.
\end{equation}
\end{definition}

\begin{thm}\label{thm2}
If an operator $G \in S_{L(X,Y)}$ is a Daugavet center then $A:=G(B_X)$ is a quasi-ball.
\end{thm}
\noindent P r o o f. 
From Theorem \ref{theo2.1}, item (iii) follows in particular
that for every $y \in Y$ and for every $\varepsilon > 0$ there is an $x
\in B_X$ with $\|y - Gx\| > \|y\| + 1 - \varepsilon.$ So, $r_y(A) > \|y\|
+ 1 - \varepsilon \ge \|y\| + r_0(A)- \varepsilon.$
\smallskip

\begin{thm}\label{thm3}
If an operator $G \in S_{L(X,Y)}$ is a Daugavet center then, for every $y \in
Y$ and for every $r \in [0, r_y(G(B_X)))$
\begin{equation} \label{eq4}
  V:= \mathop{\rm \overline{conv}}\nolimits\left( B_X \setminus G^{-1}(B_Y(y,r))\right) \supset B_X.
\end{equation}
\end{thm}
\noindent P r o o f. 
Assume it is not true. Then there are a $y \in Y$ and an $r \in
[0, r_y(G(B_X))$ such that the corresponding $V$ does not contain the
whole $B_X.$ Consider a slice $S=S(x^*,\varepsilon_0)$ of $B_X$ which
does not intersect $V$. For this slice we have $S \subset
G^{-1}(B_Y(y,r))$. Select such a small $\delta >0,$ that $\|y\| +
1 - \delta > r.$ By Theorem \ref{theo2.1}, item (iii) applied
to $x_0^*= x^*$, $y_0= -y$ and $\varepsilon = \min\{\varepsilon, \delta\}$ there
is an $x \in S \subset G^{-1}(B_Y(y,r))$ with $\|Gx - y\|
> \|y\| + 1 - \delta > r$. But in such a case $Gx \notin
B_Y(y,r),$ i.e. $x \notin G^{-1}(B_Y(y,r)),$ which leads to
contradiction.
\smallskip
\begin{cor} \label{cor1}
If an operator $G \in S_{L(X,Y)}$ is a Daugavet center then 
$A:=G(B_X)$ is anti-dentable.
\end{cor}
\noindent P r o o f. 
According to the previous theorem for every $y \in Y$ and for
every $r \in [0, r_y(A))$ the inclusion (\ref{eq4}) holds true. So
$$
A \subset G(V) \subset \mathop{\rm \overline{conv}}\nolimits G \left( B_X \setminus
G^{-1}(B_Y(y,r))\right) = \mathop{\rm \overline{conv}}\nolimits\left( A \setminus B_Y(y,r)
\right).
$$
\smallskip

Now we prove that the properties from Theorems \ref{thm2} and
\ref{thm3} together give a characterization of Daugavet centers. In fact we
prove even more:

\begin{thm}\label{theo2.3}
An operator $G \in S_{L(X,Y)}$ is a Daugavet center if and only if  it satisfies the
following two conditions:
\begin{enumerate}
\item the set $A:=G(B_X)$ is a quasi-ball; \item the condition
(\ref{eq4}) holds true for all $y \in Y$ and for all $r \in
[0, r_y(G(B_X))).$
\end{enumerate}
Moreover, if $G$ is a Daugavet center then the equation (\ref{eqDC}) holds
true for every strong Radon-Nikod\'ym operator $T \in L(X,Y)$.
\end{thm}
\noindent P r o o f. 
What remains to prove is that conditions (1) and (2) imply the
equation (\ref{eqDC}) for every strong Radon-Nikod\'ym operator $T
\in L(X,Y)$. Fix an $\varepsilon > 0.$ Let $x \in S_X$ be an element for
which $\|Tx\| > \|T\| - \varepsilon$ and $Tx$ belongs to a slice $\tilde{S}$ of $T(B_X)$ with diameter smaller than $\varepsilon$. Put $r = r_{Tx}(G(B_X)) - \varepsilon$. Then $T^{-1}\tilde{S}$ is a slice of $B_X$, so
$$
T^{-1}\tilde{S} \bigcap \left(B_X \setminus G^{-1}(B_Y(Tx,r))
\right)\neq \emptyset,
$$
Hence there is an $x_0 \in B_X$ such that $Tx_0 \in \tilde{S}$
(and, consequently, $\|Tx_0 - Tx\| < \varepsilon$) but $Tx_0 \notin
B_Y(Tx,r)$, i.e $\|Gx_0 - Tx\| > r.$ Then
$$
\|G - T\| \ge \|Gx_0 - Tx_0\| \ge \|Gx_0 - Tx\| - \varepsilon > r - \varepsilon
= r_{Tx}(G(B_X)) - 2\varepsilon
$$
$$
= \|Tx\| + r_{0}(G(B_X)) - 2\varepsilon \ge \|Tx\| + \|G\| - 2\varepsilon \ge
\|T\| + \|G\| - 3\varepsilon.
$$
\smallskip

\begin{remark} \label{rem2} 
Theorem \ref{theo2.3} will not hold true if we require $G(B_X)$ to be anti-dentable instead of the condition (2) of this Theorem. Consider $G {:}\allowbreak\  C[0,1] \oplus _1 \mathbb{R} \to C[0,1]$, $G(f, a)=f$. It is obvious that $G(B_{C[0,1] \oplus _1 \mathbb{R}})=B_{C[0,1]}$. Since $C[0,1]$ has the Daugavet property, $B_{C[0,1]}$ is an anti-dentable quasi-ball. Let us show that $G$ is not a Daugavet center. Consider rank~-~1 operator $T {:}\allowbreak\  C[0,1] \oplus _1 \mathbb{R} \to C[0,1]$, $T(f, a)=a\cdot y_0$ for some $y_0 \in S_{C[0,1]}$. So $\|G\|+\|T\|=2$, but
\begin{eqnarray*}
\|G+T\| &=& \sup_{(f,a) \in S_X}\|(G+T)(f,a)\| \\ &=& \sup_{(f,a)\in S_X}\|f+a\,y_0\| ~\leq~ \sup_{(f,a)\in S_X}(\|f\|+|a|) = 1.
\end{eqnarray*}
\end{remark}

For a set $\Gamma$ denote by ${\rm FIN}(\Gamma)$ the set of all finite
subsets of $\Gamma$. Recall that a (maybe uncountable) series
$\sum_{n\in \Gamma} x_n$ in a Banach space $X$ is said to be
unconditionally convergent to $x \in X$ if, for every $\varepsilon > 0$
there is an $A \in {\rm FIN}(\Gamma)$ such that for every $B \in
{\rm FIN}(\Gamma)$, $B \supset A$
$$
\|x - \sum_{n\in B} x_n\| < \varepsilon.
$$
\begin{thm} \label{theo2.4}
Let $G \in L(X,Y)$. Suppose that inequality $\|G+T\| \ge C+\|T\|$
with $C>0$ holds for every operator $T$ from a subspace ${\mathcal M} \subset
L(X,Y)$ of operators. Let $\widetilde{T} = \sum_{n\in \Gamma} T_n$
be a (maybe uncountable) pointwise unconditionally convergent
series of operators $T_n \in {\mathcal M}$. Then $\|G-\widetilde{T}\| \geq
C$.
\end{thm}

\noindent P r o o f. 
Pointwise unconditional convergence of $\sum_{n\in \Gamma} T_n$
implies that for every $x \in X$
$$
\sup \left\{ \biggl\| \sum_{n\in A} T_{n}x \biggr\| {:}\allowbreak\ 
                    A\in {\rm FIN}(\Gamma) \right\} < \infty.
$$
Consequently, by the Banach-Steinhaus theorem, the quantity
$$
\alpha = \sup \left\{ \biggl\| \sum_{n\in A} T_{n} \biggr\| {:}\allowbreak\ 
                    A\in {\rm FIN}(\Gamma) \right\}
$$
is finite, and whenever $B\subset \Gamma$, then
$$
\biggl\| \sum_{n\in B} T_{n} \biggr\| \le
 \sup \left\{ \biggl\| \sum_{n\in A} T_{n} \biggr\| {:}\allowbreak\ 
                    A\in {\rm FIN}(\Gamma),\ A\subset B \right\}
\le\alpha.
$$
Let $\varepsilon>0$ and pick $A_{0}\in {\rm FIN}(\Gamma)$ such that
$\|\sum_{n\in A_{0}} T_{n} \| \ge \alpha-\varepsilon$. Then we have
$$
\|G-\widetilde{T}\| \ge \biggl\| G-\sum_{n\in A_{0}} T_{n}
\biggr\| - \biggl\| \sum_{n\notin A_{0}} T_{n} \biggr\| \ge C +
\biggl\| \sum_{n\in A_{0}} T_{n} \biggr\| - \alpha \ge C-\varepsilon,
$$
which proves the theorem.
\smallskip

\begin{remark} \label{rem2.4}
Let $G {:}\allowbreak\  X \to Y$ be a non-zero Daugavet center. Since by Theorem
\ref{theo2.3} every weakly compact operator satisfies (\ref{eqDC})
the above theorem means in particular that $G$ cannot be
represented as a pointwise unconditionally convergent series of
weakly compact operators. So neither $X$ nor $Y$ can have an
unconditional basis (countable or uncountable) or can be represented
as unconditional sum of reflexive subspaces.
\end{remark}

\begin{lemma}\label{lem2.5}
Let $G {:}\allowbreak\  X \to Y$ be a Daugavet center, $\|G\|=1$. Then for every
finite-dimensional subspace $Y_{0}$ of $Y$, every $\varepsilon_{0}>0$ and
every slice $S(x^*_{0},\varepsilon_{0})$ of $B_{X}$ there is a slice
$S(x^*_{1},\varepsilon_{1}) \subset S(x^*_{0},\varepsilon_{0})$ of $B_{X}$ such
that
\begin{equation}  \label{eq2.5}
\|y+tGx\|\ge (1-\varepsilon_{0})(\|y\|+|t|) \qquad\forall y\in Y_{0},\
x\in S(x^*_{1},\varepsilon_{1}),\ \forall t \in {\mathbb R}.
\end{equation}
\end{lemma}

\noindent P r o o f. 
Let $\delta= \varepsilon_{0}/2$ and pick a finite $\delta$-net
$\{y_{1},\ldots,y_{n}\}$ in $S_{Y_{0}}$. By a repeated application
of Theorem~\ref{theo2.1}, item $(ii)$ we obtain a sequence of
slices $S(x^*_{0},\varepsilon_{0}) \supset S(u_1^*, \delta_1) \supset
\ldots \supset S(u_n^*, \delta_n)$ such that one has
\begin{equation}  \label{eq2.5.1}
\|y_{k}+Gx \| \ge 2-\delta
\end{equation}
for all $x\in S(u_k^*, \delta_k)$. Put $x^*_{1}= u_n^*$ and
$\varepsilon_{1}= \delta_n$; then (\ref{eq2.5.1}) is valid for every
$x\in S(x^*_{1},\varepsilon_{1})$ and $k=1,\ldots,n$. This implies that
for every $x\in S(x^*_{1},\varepsilon_{1})$ and every $y\in S_{Y_{0}}$
the condition
$$
\|y+Gx\| \ge 2-2\delta = 2-\varepsilon_{0}
$$
holds.

Let $0\le t_{1},t_{2}\le 1$ with $t_{1}+t_{2}=1$. If $t_{1}\ge
t_{2}$, we have for $x$ and $y$ as above
\begin{eqnarray*}
\|t_{1}Gx + t_{2}y\| &=&
\|t_{1}(Gx+y) + (t_{2}-t_{1})y\| ~\ge~
t_{1}\|Gx+y\| - |t_{2}-t_{1}|\, \|y\| \\
&\ge&
t_{1}(2-\varepsilon_{0}) + t_{2}-t_{1} ~=~ t_{1}+t_{2}-t_{1}\varepsilon_{0} ~\ge~ 1-\varepsilon_{0},
\end{eqnarray*}
and an analogous argument shows this estimate in case
$t_{1}<t_{2}$.

This implies (\ref{eq2.5}), by the homogeneity of the norm and the
symmetry of $S_{Y_{0}}$.

\smallskip

\begin{thm} \label{theo2.6}

Let $G {:}\allowbreak\  X \to Y$ be a Daugavet center. Then $G$ fixes a copy of $\ell_{1}$.

\end{thm}

\noindent P r o o f. 
Using Lemma~\ref{lem2.5} inductively, we construct sequences
of vectors $\{x_n\}_{n=1}^{\infty} \subset X$ and $\{y_n\}_{n=1}^{\infty} \subset Y$ and a sequence of slices $S(x_n^*, \varepsilon_n)$, $\varepsilon_n \leq 2^{-n}$, $n \in {\mathbb N}$, such that $y_n=Gx_n$, $x_n \in S(x_n^*, \varepsilon_n)$ and for every $y \in \mathop{\rm lin}\nolimits \{y_1, \dots , y_n \}$ and every $t \in {\mathbb R}$ the inequality
$$
\|y+ty_{n+1}\|\ge (1-\varepsilon_{n})(\|y\|+|t|\|y_{n+1}\|)
$$
holds true. Hence the sequences $\{x_n\}_{n=1}^{\infty} \subset X$ and $\{y_n\}_{n=1}^{\infty} \subset Y$ are equivalent to the canonical basis in $\ell_{1}$, and $G$ fixes a copy of $\ell_{1}$.
\smallskip

\section{Some examples of Daugavet centers}\label{sec3}

\begin{prop}\label{prop 3.1}
Let $G {:}\allowbreak\  X \to Y$ be a Daugavet center. Then for all surjective linear isometries $V {:}\allowbreak\  X \to
X$  and $U {:}\allowbreak\  Y \to Y$ the operator $UGV$ is also a Daugavet center.
\end{prop}

\noindent P r o o f. 
Let $T {:}\allowbreak\  X \to Y$ has rank one. Then 
\begin{eqnarray*}
\|UGV + T\| &=& \|U(GV + U^{-1}T)\|~=~ \| U \|\|GV + U^{-1}T \| \\ &=& \|GV + U^{-1}T \|~=~\|(G + U^{-1}TV^{-1})V \|\\ &=& \|G + U^{-1}TV^{-1}\|\|V \|~=~\|G + U^{-1}TV^{-1}\|.
\end{eqnarray*}
The operator $T$ can be represented as $Tx=x_0^*(x)y_0$ with $y_0\in Y$,
$x_0^* \in X^*$ hence $U^{-1}TV^{-1}x=x_0^*(V^{-1}x)U^{-1}y_0$ is also a rank~-~1 operator. Since G is a Daugavet center 
\begin{eqnarray*}
\|UGV + T\| &=& \|G + U^{-1}TV^{-1}\| ~=~ \|G\| + \|U^{-1}TV^{-1}\|\\ &=& \|G\| + \|T\|~=~\|UGV\| + \|T\|.
\end{eqnarray*}

\smallskip

\begin{prop}\label{prop 3.2}
Let $G {:}\allowbreak\  X \to Y$ be a Daugavet center. Then $\widetilde{G}: X/\mathop{\rm Ker}\nolimits G \to
Y$ (the natural injectivization of $G$) is also a Daugavet center.
\end{prop}

\noindent P r o o f. 
We will prove this proposition using Definition \ref{defDC} of a Daugavet center.
Let $T \in L(X/\mathop{\rm Ker}\nolimits G , Y)$ be a rank~-~1 operator and $q {:}\allowbreak\  X \to X/\mathop{\rm Ker}\nolimits G$ be the corresponding quotient mapping. Then the composition $T \circ q {:}\allowbreak\  X \to Y$ is a linear continuous rank~-~1 operator. Since $G$ is a Daugavet center the identity $$\|G + T \circ q\| = \|G\| + \|T \circ q\|$$ holds true. The operator $\widetilde{G}$ is the natural injectivization of $G$ hence $\|G\|= \|\widetilde{G}\|$. It is well-known that $q(\stackrel{_\circ}{B}_{X})=\stackrel{_\circ}{B}_{X/\mathop{\rm Ker}\nolimits G}$ where $\stackrel{_\circ}{B}_{X}$ and $\stackrel{_\circ}{B}_{X/\mathop{\rm Ker}\nolimits G}$ are the open unit balls of $X$ and $X/\mathop{\rm Ker}\nolimits G$ respectively. This implies that $\|T \circ q\|=\|T\|$ and $\|\widetilde{G} + T \|=\|(\widetilde{G} + T)\circ q \|=\|G + T \circ q\|$. So we have $$\|\widetilde{G} + T \| =\|G + T \circ q\| = \|G\| + \|T \circ q\|=\|\widetilde{G}\| + \|T \|$$
which proves the proposition.
\smallskip

\begin{lemma}\label{lem3.3}
If $G_1 {:}\allowbreak\  X_{1} \to Y_{1}$ and $G_2 {:}\allowbreak\  X_{2} \to Y_{2}$ are Daugavet centers, $\|G_1\|=\|G_2\|=1$.
Then operator $G {:}\allowbreak\  X_{1}\oplus_{\infty}X_{2} \to Y_{1}\oplus_{\infty}Y_{2}$ ($G {:}\allowbreak\  X_{1}\oplus_{1}X_{2} \to Y_{1}\oplus_{1}Y_{2}$) which maps every $(x_1, x_2)$ into $(G_1x_1, G_2x_2)$ is a Daugavet center.
\end{lemma}

\noindent P r o o f. 
We first prove that $G {:}\allowbreak\  X_{1}\oplus_{\infty}X_{2} \to Y_{1}\oplus_{\infty}Y_{2}$ is a Daugavet center.
Consider
$x^*_{j}\in X_{j}^{*}$, $y_{j}\in Y_{j}$ ($j=1,2$) with
$\|(y_{1},y_{2})\|=\max\{\|y_{1}\|,\|y_{2}\|\} =1$,
$\|(x^*_{1},x^*_{2})\| = \|x^*_{1}\| + \|x^*_{2}\| =1$.
Assume without loss of generality that $\|y_{1}\|=1$. We will use the characterization of Daugavet centers from item (iii) of 
Theorem~\ref{theo2.1}.
For a given $\varepsilon>0$ there is an $x_{1}\in X_{1}$
satisfying
$$
\|x_{1}\|=1,\quad x^*_{1}(x_{1})\ge \|x^*_{1}\|(1-\varepsilon) , \quad
\|G_{1}x_{1}+y_{1}\|\ge 2-\varepsilon.
$$
Also, pick $x_{2}\in X_{2}$ such that
$$
\|x_{2}\|=1,\quad x^*_{2}(x_{2})\ge \|x^*_{2}\|(1-\varepsilon).
$$
Then $\|(x_{1},x_{2})\|=1$, $\langle (x^*_{1},x^*_{2}), (x_{1},x_{2})
\rangle \ge 1-\varepsilon$ and
$$
\|G(x_{1},x_{2}) + (y_{1},y_{2}) \| \ge
\|Gx_{1} + y_{1}\| \ge 2-\varepsilon.
$$
Thus, $G$ is a Daugavet center.

A similar calculation, based on item (v) of Theorem~\ref{theo2.1}, 
proves that $G {:}\allowbreak\  X_{1}\oplus_{1}X_{2} \to Y_{1}\oplus_{1}Y_{2}$ is a Daugavet center.
\smallskip

\begin{lemma}\label{lem3.4}
Let $G {:}\allowbreak\  X \to Y$ be a Daugavet center, $\|G\|=1$. Denote $\widetilde{G} {:}\allowbreak\  X \to Y\oplus_{1}Y_{1}$, $\widetilde{G}x=(Gx, 0)$, and $\hat{G} {:}\allowbreak\  X_{1} \oplus_{\infty}X \to Y$, $\hat{G}(x_1, x_2)=Gx_{2}$. Then:
\begin{statements}
  \item the operator $\widetilde{G}$ is a Daugavet center;
  \item the operator $\hat{G}$ is a Daugavet center.
\end{statements}
\end{lemma}

\noindent P r o o f. 
Part (b) can be proved in a similar manner as Lemma~\ref{lem3.3}, so we only present the proof of (a).
Consider
$x^*\in S_{X^{*}}$, $y_{j}\in Y_{j}$ ($j=0,1$) with
$\|(y_{0},y_{1})\|=\|y_{0}\|+\|y_{1}\| =1$. By
Theorem~\ref{theo2.1}
there is, given $\varepsilon>0$, some $x_{0}\in S(x^*, \varepsilon)$
satisfying 
$$
\biggl\|Gx_{0}+\frac{y_{0}}{\|y_{0}\|}\biggr\|\ge 2-\varepsilon.
$$
Then we have
\begin{eqnarray*}
\|\widetilde{G}x_{0}+(y_{0},y_{1})\| &=& \|Gx_{0} + y_{0}\| + \|y_{1}\| \\
&\geq & \biggl\|Gx_{0} + \frac{y_{0}}{\|y_{0}\|} + y_{0}\biggl(1 - \frac{1}{\|y_{0}\|}\biggr)\biggr\| + \|y_{1}\| \\
&\geq & \biggl\|Gx_{0} + \frac{y_{0}}{\|y_{0}\|}\biggr\| + \|y_{0}\|\biggl(1 - \frac{1}{\|y_{0}\|}\biggr) + \|y_1\| \\
&\geq & 2- \varepsilon
\end{eqnarray*}
which proves the Lemma.
\bigskip

Let $K$ be a compact space without isolated points. Then $C(K)$ has the Daugavet property and this means that the identity operator is a Daugavet center. Therefore by Proposition \ref{prop 3.1} every surjective linear isometry $V {:}\allowbreak\  C(K) \to C(K)$ is a Daugavet center.

In particular, if we consider any bijective continuous function $\varphi {:}\allowbreak\  K \to K$ then the operator $G_{\varphi} {:}\allowbreak\ C(K) \to C(K)$, $G_{\varphi}f = f \circ \varphi $ is a surjective linear isometry and hence a Daugavet center.

Our next aim is to prove that for every continuous function $\varphi {:}\allowbreak\  K \to K$ such that $\varphi^{-1}(t)$ is nowhere dense in $K$ for all $t \in K$ the corresponding operator $G_{\varphi}$ is a Daugavet center as well.

\begin{lemma}\label{lem3.5}
For an operator $G {:}\allowbreak\  X \to C(K)$, $\|G \|=1$, the following
assertions are equivalent:
\begin{aequivalenz}
\item $G$ is a Daugavet center.
\item For every $\varepsilon >0$, every open set $U \subset K$, every $x^* \in S_{X^*}$ and $s= \pm 1$ there is $f \in S(x^*, \varepsilon)$ such that 
$$ \sup_{t \in U} s \cdot (Gf)(t) > 1-\varepsilon . $$
\end{aequivalenz}
\end{lemma}

\noindent P r o o f. 
$(i)\Rightarrow(ii)$ Let us consider a function $g \in S_{C(K)}$ such that $\mathop{\rm supp}\nolimits g \subset U$ and $s \cdot g \geq 0$. By Theorem~\ref{theo2.1}
for every $\varepsilon > 0$ and every $x^* \in S_{X^*}$ there is an element $f \in S(x^*, \varepsilon )$ such that 
$$\sup_{t \in K} |(Gf+g)(t)| > 2- \varepsilon .$$

Remark that   
$|Gf+g|= |Gf| \leq 1$ on $K \setminus U$
and hence $|Gf+g|$ attains its supremum in 
$U$. Then there is a point $t_0 \in U$ which fulfills the inequality $|(Gf+g)(t_0)| > 2-\varepsilon$. Since $s \cdot g \geq 0$, then $s \cdot (Gf)(t_0) \geq 0$ and  
\begin{eqnarray*}
|(Gf+g)(t_0)| &=& s \cdot (Gf+g)(t_0) ~\leq~ \sup_{t \in U} s \cdot (Gf+g)(t) \\
 &\leq & \sup_{t \in U} s \cdot (Gf)(t) + \sup_{t \in U} s \cdot g(t) ~\leq~ \sup_{t \in U} s \cdot (Gf)(t) +1.
\end{eqnarray*}
Therefore $$ \sup_{t \in U} s \cdot (Gf)(t) > 1-\varepsilon . $$
$(ii)\Rightarrow(i)$ Consider $g \in S_{C(K)}$, pick some $\tau \in K$ with $|g(\tau)|=1$ and put $s = g(\tau)$. Then for every $\varepsilon > 0$ there is an open neighborhood $U$ of $\tau$ such that $s \cdot g > 1- \varepsilon $ on $U$.
For every $x^* \in S_{X^*}$ there is an element $f \in S(x^*, \varepsilon)$ which satisfies the inequality 
$$ \sup_{t \in U} s \cdot (Gf)(t) > 1-\varepsilon . $$
Then we have 
\begin{eqnarray*}
\|Gf+g\| &=& \sup_{t \in K} |s \cdot (Gf+g)(t)|~\geq~ \sup_{t \in U} (s \cdot (Gf)(t)+ s \cdot g(t)) \\ & \geq & \sup_{t \in U} s \cdot (Gf)(t)+ 1 - \varepsilon ~\geq~ 1 - \varepsilon +1 - \varepsilon = 2 - 2\varepsilon .
\end{eqnarray*}
By Theorem~\ref{theo2.1} 
$G$ is a Daugavet center.
\bigskip

Consider in Lemma \ref{lem3.5} $X=C(K_{1})$. By the Riesz representation theorem for any linear functional $x^*$ on $C(K_{1})$, there is a unique Borel regular signed measure $\sigma$ on $K_{1}$ such that 
$$x^*(f) = \int\limits_{K_{1}} f\, d\sigma $$
for all $f \in C(K_{1})$, and $\|x^*\| = |\sigma|(K_{1})$. So every slice
\begin{eqnarray*} 
S(x^*, \varepsilon) &=& \{ f \in B_{C(K_{1})} {:}\allowbreak\  \int\limits_{K_{1}} f\, d\sigma \geq | \sigma |(K_{1}) - \varepsilon \} \\
&=& \{ f \in B_{C(K_{1})} {:}\allowbreak\  \int\limits_{K_{1}}(1-f\, ({\bf{1}}_{K_{1}^+}-{\bf{1}}_{K_{1}^-}))\, d|\sigma | \leq \varepsilon \}. 
\end{eqnarray*}
Here $K_{1}=K_{1}^+ \sqcup K_{1}^-$ is a Hahn decomposition of $K_{1}$ for $\sigma $ and ${\bf{1}}_A $ denotes a characteristic function of the set $A$.

So in the case of $X=C(K_{1})$ Lemma \ref{lem3.5} reformulates as follows:

\begin{lemma} \label{lem3.9}
For an operator $G {:}\allowbreak\  C(K_{1}) \to C(K_{2})$, $\|G\|=1$, the following
assertions are equivalent:
\begin{aequivalenz}
\item $G$ is a Daugavet center.
\item For every $\varepsilon >0$, every open set $U \subset K_{2}$ and every Borel regular signed measure $\sigma$ on $K_{1}$
and $s= \pm 1$ 
there is an $f \in B_{C(K_{1})}$ such that
\begin{equation} \label{eq3.1}
\int\limits_{K_{1}}(1-f\, ({\bf{1}}_{K_{1}^+}-{\bf{1}}_{K_{1}^-}))\, d|\sigma | \leq \varepsilon
\end{equation}
and 
\begin{equation} \label{eq3.2}
\sup_{t \in U} s \cdot (Gf)(t) > 1-\varepsilon .
\end{equation}
\end{aequivalenz}
\end{lemma} 

\begin{thm} \label{theo3.6}
Let $K_{1}$ and $K_{2}$ be compact spaces without isolated points, $\varphi {:}\allowbreak\  K_{2} \to K_{1}$ be a continuous function such that for every $t \in K_{1}$ the set $\varphi^{-1}(t)$ is nowhere dense in $K_{2}$. Suppose that an operator $G_{\varphi} {:}\allowbreak\  C(K_{1}) \to C(K_{2})$ maps every $f \in C(K_{1})$ into the composition $f \circ \varphi $. Then $G_{\varphi}$ is a Daugavet center.
\end{thm}

\noindent P r o o f. 
Consider an $\varepsilon >0$, an open set $U \subset K_{2}$, and a Borel regular signed measure $\sigma$ on $K_{1}$, and put $s=1$. We will construct a function $f \in B_{C(K_{1})}$ satisfying (\ref{eq3.1}) and (\ref{eq3.2}). 

The measure $\sigma$ can have at most countable set of atoms. Let us show that for every open  $U \subset K_{2}$ the set $\varphi (U)$ is uncountable. Assume that there exists an open set $U \subset K_{2}$ for which it is not true. Then $\varphi^{-1} (\varphi (U))$ is a countable union of nowhere dense sets in $K_{2}$ because for every $t \in \varphi (U) \subset K_{1}$ the set $\varphi^{-1} (t)$ is nowhere dense in $K_{2}$ by the condition of this Theorem. This contradicts the Baire category theorem.

So we can pick a point $t_0 \in U$ such that $\varphi (t_0)$ is not an atom of  $\sigma $,  i. e. $|\sigma |(\varphi (t_0))=0$. Moreover, since $\sigma $ is a Borel regular measure, there is an open neighborhood $V \subset \varphi (U)$ of the point $\varphi (t_0)$ such that $|\sigma|(V) < \varepsilon / 4$.

Now we pass on to the construction of $f$. To satisfy (\ref{eq3.2}) we select $f$ in such a way that $f(\varphi(t_0)) > 1-\varepsilon$. First we pick a function $\tilde{f} \in S(\sigma, \varepsilon /2)$. If $\tilde{f}(\varphi (t_0))>1- \varepsilon$, then we can simply put $f = \tilde{f}$.

If $\tilde{f}(\varphi (t_0)) \leq 1- \varepsilon$, we put $f=\tilde{f}$ in 
\hbox{$K_{1} \setminus V$} and \hbox{$f(\varphi (t_0))=1$}. Since $K_{1} \setminus V \cup {\{ \varphi (t_0) \}}$ is closed, we can use the Tietze extension theorem to construct a continuous extension  $f$ on 
$V \setminus \varphi (t_0)$ and keep the condition $\|f\|=1$. Now we show that (\ref{eq3.1}) also holds for this $f$:
\begin{eqnarray*}
\int\limits_{K_{1}}(1-f\, ({\bf{1}}_{K_{1}^+} - {\bf{1}}_{K_{1}^-} ))\, d|\sigma|&=&
\int\limits_{K_{1}\setminus V}(1-\tilde{f}\, ( {\bf{1}}_{K_{1}^+} - {\bf{1}}_{K_{1}^-} ))\, d|\sigma| \\
&+&\int\limits_{V}(1-f \, ({\bf{1}}_{K_{1}^+} -{\bf{1}}_{K_{1}^-} ))\, d|\sigma| \\
&\leq&
\int\limits_{K_{1}}(1-\tilde{f}\, ({\bf{1}}_{K_{1}^+} -{\bf{1}}_{K_{1}^-} ))\, d|\sigma| + 
\varepsilon /2 \\
&\leq& \varepsilon /2 + \varepsilon /2 ~=~ \varepsilon .
\end{eqnarray*}

So for every $\varepsilon >0$, every Borel regular measure $\sigma$ on $K_{1}$, every open set $U \subset K_{2}$ and $s=1$ we have a function $f \in B_{C(K_{1})}$ satisfying the inequalities (\ref{eq3.1}) and (\ref{eq3.2}). The case $s=-1$ can be proved in a very similar way. Thus, by Lemma \ref{lem3.9} $G_{\varphi}$ is a Daugavet center. 
\bigskip

Let us give an example of a Daugavet center on $C(K)$ of a very different nature.

\begin{prop} \label{prop3.8}
Consider $K_{1}=[-1; 1]$ and define $G {:}\allowbreak\  C(K_{1}) \to C(K_{1})$ as $(Gf)(x) = \frac {f(x)+f(-x)}{2}$. Then $G$ is a Daugavet center.
\end{prop}

\noindent P r o o f. 
We will use Lemma \ref{lem3.9} to prove this proposition. Let us fix an $\varepsilon > 0$, a Borel regular signed measure $\sigma$, an open set $U \subset K_{1}$, $s=1$ and a function $\tilde{f} \in S(\sigma , \varepsilon /2)$. If there is a point $t_0 \in U $ such that $\frac{\tilde{f}(t_0)+\tilde{f}(-t_0)}{2} > 1-\varepsilon $ then
$$ \sup_{t\in U} s\cdot (G\tilde{f})(t) > 1 - \varepsilon .$$
Otherwise we pick a point $t_1 \in U$ such that neither $t_1$ nor $-t_1$ is an atom of $\sigma$. Consider disjoint segments $[a_1, b_1], [a_2, b_2] \subset K_{1}$ 
such that $|\sigma|([a_1, b_1])<\varepsilon/8$, $t_1 \in [a_1, b_1]$ and $|\sigma|([a_2, b_2])<\varepsilon/8$, $-t_1 \in [a_2, b_2]$. Let $\tilde{f}_1 {:}\allowbreak\  [a_1, b_1] \to K_{1}$ be a continuous function such that $\tilde{f}_1(a_1)=\tilde{f}(a_1)$, $\tilde{f}_1(b_1)=\tilde{f}(b_1)$ and $\tilde{f}_1(t_1)=1$. Let $\tilde{f}_2 {:}\allowbreak\  [a_2, b_2] \to K_{1}$ be a continuous function such that $\tilde{f}_2(a_2)=\tilde{f}(a_2)$, $\tilde{f}_2(b_2)=\tilde{f}(b_2)$ and $\tilde{f}_2(-t_1)=1$. Then denote $\Delta := K_{1} \setminus \{[a_1, b_1] \cup [a_2, b_2]\}$ put 
$$f:={\bf{1}}_{[a_1, b_1]}\tilde{f}_1+{\bf{1}}_{[a_2, b_2]}\tilde{f}_2+{\bf{1}}_{\Delta}\tilde{f}.$$
Then $f \in B_{C(K_{1})}$ and 
\begin{eqnarray*}
&&\int\limits_{K_{1}}(1-f\, ( {\bf{1}}_{K_{1}^+} - {\bf{1}}_{K_{1}^-} ))\, d|\sigma| ~=~  \int\limits_{\Delta}(1-\tilde{f}\, ( {\bf{1}}_{K_{1}^+} - {\bf{1}}_{K_{1}^-} ))\, d|\sigma| \\&+& \int\limits_{[a_1, b_1]}(1-\tilde{f_1}\, ( {\bf{1}}_{K_{1}^+} - {\bf{1}}_{K_{1}^-} ))\, d|\sigma|~+~ \int\limits_{[a_2, b_2]}(1-\tilde{f_2}\, ( {\bf{1}}_{K_{1}^+} - {\bf{1}}_{K_{1}^-} ))\, d|\sigma| \\ &\leq& \int\limits_{K_{1}}(1-\tilde{f}\, ( {\bf{1}}_{K_{1}^+} - {\bf{1}}_{K_{1}^-} ))\, d|\sigma| ~+~ \varepsilon /4 ~+~ \varepsilon /4 ~\leq~  \varepsilon /2 + \varepsilon /2  ~=~ \varepsilon .   
\end{eqnarray*}
Hence $f \in S(\sigma, \varepsilon)$ and $$\sup_{t\in U} s\cdot (Gf)(t)\geq \frac{f(t_1)+f(-t_1)}{2} =1.$$
If we put $s=-1$, the analogous conclusions prove the proposition.
\bigskip

A quite non-trivial class of Daugavet centers was discovered in \cite{Pop}: every isometric embedding $G {:}\allowbreak\  L_{1}[0,1] \to L_{1}[0,1]$ is a Daugavet center. Let us show that the analogous result for $C[0,1]$ is false. This will answer in negative a question from \cite{Pop1}.

\begin{example} \label{exp3}
Consider $T {:}\allowbreak\  C[0,1] \to C[0,1]$,
$$
Tf=
\left\{
\begin{array}{r}
f(2t)\mbox{\quad if \quad}t \in \left[0, \frac{1}{2}\right], \\
2f(1)(1-t)\mbox{\quad if \quad}t \in \left(\frac{1}{2}, 1\right].
\end{array}
\right.
$$
Let us prove that $T$ is an isometric embedding. It is obvious that $T$ is a linear operator. Remark that $|Tf|$ attains its supremum in $\left[0, \frac{1}{2}\right]$ because for every $t \in \left(\frac{1}{2}, 1\right]$ we have $|Tf(t)|= |2f(1)(1-t)| < |f(1)|=|Tf\left( \frac{1}{2} \right)|$. Hence for every $f \in C[0,1]$ $$\|Tf\|=\sup_{t \in \left[0, \frac{1}{2}\right]} |Tf(t)| = \sup_{t \in [0,1]} |f(t)|= \|f\|.$$

Now we show with the help of Lemma \ref{lem3.5} that $T$ is not a Daugavet center. Our aim is to find an $\varepsilon>0$, an open set $U \subset [0,1]$ and an $x^* \in S_{C^*[0,1]}$ such that every $f \in S(x^*, \varepsilon)$ satisfies $\sup_{t \in U} Tf(t) \leq 1- \varepsilon$. If we put $\varepsilon := \frac{1}{4}$ and $U := \left(\frac{3}{4}, 1\right]$ then for every $f \in B_{C[0,1]}$ we have $$\sup_{t \in U} Tf(t) = \sup_{t \in \left(\frac{3}{4}, 1\right]} 2f(1)(1-t) \leq 2|f(1)|\left( 1-\frac{3}{4} \right) = \frac{|f(1)|}{2} \leq \frac{1}{2} < 1- \varepsilon.$$
\end{example} 

\smallskip


\section{The main result} \label{main}
\begin{definition}\label{defmain1}
Let $E$ be a seminormed space, $A \subset B_E$, $\mathcal{U}$ be a free
ultrafilter on a set $\Gamma$, and $f: \Gamma \to A$ be a
function. The triple $(\Gamma, \mathcal{U}, f)$ is said to be an
$A$-valued $E$-atom if for every $w \in E$
\begin{equation} \label{eqmain1}
\lim_{\mathcal{U}}\|f + w\| = 1 + \|w\|.
\end{equation}
\end{definition}

The following characterization of Daugavet centers is a consequence of
Theorem~\ref{theo2.1} and Lemma \ref{lem2.5}.

\begin{thm}\label{maintheo-1}
Let $X, Y$ be Banach spaces. An operator $G \in S_{L(X,Y)}$ is a Daugavet center
if and only if  for every slice $S$ of $B_X$ there is a $G(S)$-valued
$Y$-atom.
\end{thm}

\noindent P r o o f. 
Let us start with ``if" part. We are going to prove that $G$
satisfies condition (iii) of Theorem~\ref{theo2.1}. Fix $y_0\in
S_Y$, $x_0^* \in S_{X^*}$ and $\varepsilon >0$. Denote $S=S(x_0^*,
\varepsilon)$. Due to our assumption there is a $G(S)$-valued $Y$-atom
$(\Gamma, \mathcal{U}, f)$. Plugging $w=y_0$ in (\ref{eqmain1}) we get in
particular that $\|f(t) + y_0\| > 2 - \varepsilon$ for some $t \in
\Gamma$. Since $f(t) \in G(S)$, there is an $x \in S$ such that
$f(t)= Gx$. This $x$ fulfills the required conditions $x_0^*(x)
\geq 1-\varepsilon$ and $\|Gx+y_0\|>2-\varepsilon$. The ``if" part is proved.

Let us demonstrate the ``only if" part. Fix a slice $S$ of $B_X$.
Put $\Gamma = {\rm FIN}(Y)$ and take the natural filter $\mathcal{F}$ on
$\Gamma$ whose base is formed by the collection of subsets $\hat A
\subset {\rm FIN}(Y)$, $A \in {\rm FIN}(Y)$, where $\hat A:= \{B \in
{\rm FIN}(Y){:}\allowbreak\  A \subset B\}$. According to Lemma \ref{lem2.5} for
every $A \in {\rm FIN}(Y)$ there is an element $x(A) \in S$ such that
for all $y \in A$
$$
\|y + G(x(A))\| > \left(1-\frac{1}{|A|}\right)(\|y\|+ 1).
$$
Define $f(A):=G(x(A))$. It is easy to see that for every
ultrafilter $\mathcal{U} \succ \mathcal{F}$ the triple $(\Gamma, \mathcal{U}, f)$ is the
required $G(S)$-valued $Y$-atom.
\bigskip

It is clear, that if $A \subset B$, then every $A$-valued $E$-atom
is at the same time a $B$-valued $E$-atom. A  $B_E$-valued
$E$-atom will be called just an $E$-atom.

\begin{lemma}\label{lem08-1}
Let $(E,p)$ be a seminormed space, $Y$ be a subspace of $E$ and
$(\Gamma, \mathcal{U}, f)$ be a $Y$-atom. Define
$$
p_r(x)= \mathcal{U}\mbox{-}\lim_{t} p(x + r f(t))-r
$$
for $x \in E$ and $r>0$. Then:
\begin{statements}
  \item $0\le p_r(x)\le p(x)$
for all $x\in E$,
 \item $p_r(y)=p(y)$ for all $y \in
Y$,
 \item $x \mapsto p_r(x)$ is convex for each
$r$,
 \item $r \mapsto p_r(x)$ is convex for each $x$,
  \item
$p_r(tx)= tp_{r/t}(x)$ for each $t>0$.
\end{statements}
\end{lemma}

\noindent P r o o f. 
The only thing that is not obvious is that $p_r \geq 0$; note
that (b) is just the definition of $Y$-atom. Now, given $\varepsilon>0$,
pick $t_\varepsilon$ such that $p(f(t_\varepsilon)) > 1 - \varepsilon$ and
$$
p(x+rf(t_\varepsilon)) \le \mathcal{U}\mbox{-}\lim_t p(x+r f(t)) + \varepsilon.
$$
Then
\begin{eqnarray*}
\mathcal{U}\mbox{-}\lim_t p(x+r f(t)) &\ge&
\mathcal{U}\mbox{-}\lim_t p(-r f(t_\varepsilon) + r f(t)) -  p(x+rf(t_\varepsilon)) \\
&=&
r p(f(t_\varepsilon)) + r - p(x+rf(t_\varepsilon)) \\
& \ge& 2r - r \varepsilon - \mathcal{U}\mbox{-}\lim_t p(x+r f(t)) - \varepsilon;
\end{eqnarray*}
hence $\mathcal{U}\mbox{-}\lim_t p(x+r f(t)) \ge  \frac12(2r- \varepsilon
- r \varepsilon)$ and $p_r(x) \ge 0$.
\smallskip

\begin{lemma}\label{lem08-2}
Assume the conditions of Lemma~\ref{lem08-1}.
Then $r\mapsto p_r(x)$ is decreasing  for each $x$. The quantity
$$
\bar p (x):= \lim_{r\to\infty}  p_r(x) = \inf_r p_r(x)
$$
satisfies \textrm{(a)--(c)} of Lemma~\ref{lem08-1} and moreover
\begin{equation}\label{eq08-06}
  \bar p(tx) = t \bar p(x) \qquad\mbox{for }t>0,\ x\in X.
\end{equation}
\end{lemma}

\noindent P r o o f. 
By Lemma~\ref{lem08-1}(a) and (d),
$r\mapsto p_r(x)$ is bounded and convex, hence
decreasing. Therefore, $\bar p$ is well defined. Clearly,
(\ref{eq08-06}) follows from (e) above.
\bigskip


\textbf{Proof of the main theorem} (Theorem \ref{theo-intr-2}).
Let $G {:}\allowbreak\  X \to Y$ is a Daugavet center, $Y$ be a subspace of a Banach space \ $E$,
and $J: Y \to E$ is the natural embedding operator.

Let $\mathcal{P}$ be the family of all seminorms $q$ on $E$ that
are dominated by the norm of $E$ and for which $q(y)= \|y\|$ for
$y\in Y$. By Zorn's lemma, $\mathcal{P}$ contains a minimal
element, say~$p$.

 
\begin{clm}
Every $Y$-atom $(\Gamma, \mathcal{U}, f)$ is at the same
time an $(E, p)$-atom, i.e. for every $w \in E$
\begin{equation} \label{eqmain3}
\lim_{\mathcal{U}}p(f + w) = 1 + p(w).
\end{equation}
\end{clm}

\noindent P r o o f.
To prove the claim associate to $p$ and $(\Gamma, \mathcal{U}, f)$ the
functional $\bar p$ from Lemma~\ref{lem08-2}. Note that $0\le \bar
p\le p$, but $\bar p$ need not be a seminorm. However,
$$
q(x) = \frac{ \bar p(x) + \bar p(-x)}2
$$
defines a seminorm, and $q\le p$. By Lemma~\ref{lem08-1}(b) and by
minimality of $p$ we get that
\begin{equation}\label{eq08-09}
  q(x)=p(x) \qquad \forall x\in X.
\end{equation}
Now, since $p(x)\ge \bar p(x)$ and $p(x)=p(-x)\ge \bar p(-x)$,
(\ref{eq08-09}) implies that $p(x)=\bar p(x)$. Finally, by
Lemma~\ref{lem08-1}(a) and the definition of $\bar p$ we have
$p(x)=p_r(x)$ for all $r>0$; in particular $p(x)=p_1(x)$, which is
our claim~(\ref{eqmain3}). 
\bigskip

Now let us introduce a new norm on $E$ as
$$
|\mkern-2mu|\mkern-2mu| x |\mkern-2mu|\mkern-2mu| := p(x) + \|[x]\|_{E/Y};
$$
and let us show that this is the equivalent norm that we need.
Indeed, clearly $|\mkern-2mu|\mkern-2mu| x|\mkern-2mu|\mkern-2mu| \le 2\|x\|$. On the other hand,
$|\mkern-2mu|\mkern-2mu| x|\mkern-2mu|\mkern-2mu| \ge \frac13 \|x\|$. To see this assume $\|x\|=1$. If
$\|[x]\|_{E/Y}\ge \frac13$, there is nothing to prove. If not,
pick $y\in Y$ such that $\|x-y\|<\frac13$. Then
$p(y)=\|y\|>\frac23$, and
$$
|\mkern-2mu|\mkern-2mu| x |\mkern-2mu|\mkern-2mu| \ge p(x) \ge p(y)- p(x-y) > \frac23 - \|x-y\| >
\frac13.
$$
Therefore, $\| \cdot \|$ and $|\mkern-2mu|\mkern-2mu| \cdot |\mkern-2mu|\mkern-2mu|$ are equivalent norms.
Also evidently for $y \in Y$
 $$
|\mkern-2mu|\mkern-2mu| y |\mkern-2mu|\mkern-2mu| = p(y) = \|y\|.
 $$
 What remains to prove is that $J \circ G :
X \to E$ is a Daugavet center. This can be done easily with the help of
Theorem \ref{maintheo-1} and of the Claim. Namely, let $S$ be an
arbitrary slice of $B_X$. Due to Theorem \ref{maintheo-1} it
is sufficient to demonstrate the existence of a $G(S)$-valued $(E,
|\mkern-2mu|\mkern-2mu| \cdot |\mkern-2mu|\mkern-2mu|)$-atom. Since $G : X \to Y$ is a Daugavet center, the same
Theorem \ref{maintheo-1} ensures the existence of a $G(S)$-valued
$Y$-atom $(\Gamma, \mathcal{U}, f)$. But according to the Claim, $(\Gamma,
\mathcal{U}, f)$ is also an $(E, p)$-atom. Consequently for every $w \in E$
\begin{eqnarray*}
 \lim_{\mathcal{U}} |\mkern-2mu|\mkern-2mu| f + w |\mkern-2mu|\mkern-2mu| &=& \lim_{\mathcal{U}}\left(p(f + w) + \|[f + w]\|_{E/Y}\right)
\\
 &=& \lim_{\mathcal{U}}p(f + w) + \|[w]\|_{E/Y} = 1 + p(w) + \|[w]\|_{E/Y} = 1 + |\mkern-2mu|\mkern-2mu| w |\mkern-2mu|\mkern-2mu|.
\end{eqnarray*}
This means that $(\Gamma, \mathcal{U}, f)$ is the required $G(S)$-valued
$(E, |\mkern-2mu|\mkern-2mu| \cdot |\mkern-2mu|\mkern-2mu|)$-atom. \textbf{The main theorem is
proved}.
The same renorming idea is applicable to the theory of $\ell_1$-types \cite{KSheWer}.
\bigskip

The next corollary improves the statement of remark \ref{rem2.4}.

\begin{cor} \label{remmain1}
If $G {:}\allowbreak\  X \to Y$ is a non-zero Daugavet center, then neither $X$ nor $Y$
can be embedded into a space $E$, in which the identity operator
$\mathrm{Id}_E$ has a representation as a pointwise unconditionally
convergent series of weakly compact operators. In particular
neither $X$ nor $Y$ can be embedded into a space $E$ having an
unconditional basis (countable or uncountable) or having a
representation as unconditional sum of reflexive subspaces.
\end{cor}

\noindent P r o o f.  Let $\mathrm{Id}_E = \sum_{n\in \Gamma} T_n$, where the series
is pointwise unconditionally convergent, and all the $T_n {:}\allowbreak\  E
\to E$ are weakly compact. At first assume $Y \subset E$, and
denote $J \in L(Y,E)$ the natural embedding operator. Equip $E$
with the equivalent norm from Theorem \ref{theo-intr-2} making
$J \circ G$ a Daugavet center. Then $J \circ G = \sum_{n\in \Gamma} T_n \circ
J \circ G$, the series is pointwise unconditionally convergent,
and all the operators $T_n \circ J \circ G$ are weakly compact.
This contradicts Theorem \ref{theo2.4}.

Now assume $X \subset E$. Recall that for a set $\Delta$ of big
cardinality (say, for $\Delta = B_{Y^*}$), there is an isometric
embedding $J {:}\allowbreak\  Y \to \ell_\infty(\Delta)$. Since
$\ell_\infty(\Delta)$ is an injective space (i.e. the Hahn-Banach
extension theorem holds true for $\ell_\infty(\Delta)$-valued
operators), there is an operator $U {:}\allowbreak\  E \to
\ell_\infty(\Delta)$ such that $U|_X = J \circ G$. Then
$$
J \circ G = (U \circ \mathrm{Id}_E)|_X = \sum_{n\in \Gamma} U \circ
(T_n)|_X.
$$
This representation leads to contradiction the same way as in the
previous case.
\bigskip


\end{document}